\documentclass[12pt]{amsart}

\newtheorem{theorem}{Theorem} 
\newtheorem{lemma}{Lemma}[section]

\theoremstyle{definition} 
\newtheorem{definition}[lemma]{Definition}

\newcommand{\inv}{^{-1}} 
\newcommand{\edge}[1]{\buildrel #1
\over\rightarrow } 
\newcommand{\abs}[1]{\vert #1\vert}
\newcommand{\ovr}[1]{\overline{#1 } } 
\newcommand{\hatt}[1]{\widehat{ #1 } } 
\newcommand{\set}[1]{\{#1 \}} 
\newcommand{\subgroup}[1]{\langle #1 \rangle}
 
\newcommand{\A}{\mathcal A}
\newcommand{\B}{\mathcal B} 
\newcommand{\C}{\mathcal C}
\newcommand{\e}{\epsilon}       % The empty word
\newcommand{\G}{\Gamma} 

\newcommand{\Si}{\Sigma}
\newcommand{\SSi}{{\Sigma^*}}

\title{Automatic Quotients of Free Groups} \author{Robert H. Gilman}

\address{Department of Mathematical Sciences, Stevens Institute of
Technology, Hoboken, New Jersey} \email{rgilman@math.stevens.edu}

\thanks{The author thanks Imperial College for its hospitality while
this paper was being written.}



\begin{document} 

\begin{abstract}
Automatic groups admitting prefix closed automatic structures with
uniqueness are characterized as the quotients of free groups by normal
subgroups possessing sets of free generators satisfying certain
language-theoretic conditions.
\end{abstract}

\maketitle

\section{Introduction} %%%%%%%%%%%%%%%%%%%%%%%%%%%%%%%%%%%%%%%%%%%%

It is an open problem whether or not every synchronously automatic
group admits a prefix closed automatic structure with
uniqueness~\cite[Open Problem~2.5.10]{E+}. Motivated by this problem,
we give new characterizations of synchronous and asynchronous
automatic groups which have prefix closed automatic structures with
uniqueness.

Automatic groups are a class of finitely presented groups modeled on
the fundamental groups of compact 3--manifolds and defined by the
property that group multiplication can be carried out by finite
automata.  The standard introduction is~\cite{E+}; \cite{CC} is a more
recent introduction for those familiar with the theory of finite
automata. Accounts of the context in which automatic groups occur are
given in~\cite{JWC} and~\cite{KO}.

\begin{theorem}\label{asynchronous} A group has a prefix-closed
  asynchronous automatic structure with uniqueness if and only if it is
  isomorphic to the quotient of a finitely generated free group $F$ by a
  normal subgroup $N$ admitting a linear language $L$ of freely
  reduced generators with significant letters.
\end{theorem}

$L$ is a subset of the free monoid over a set of free
generators and their inverses for $F$. Linear languages lie
between the better known classes of regular and context free
languages. $L$ has significant letters if each word in $L$ has a
distinguished letter such that free reduction of a product of two
words from $L$ or their inverses does not affect the distinguished letters
except when the product consists of a word times its inverse. Precise
definitions are given in the next section.

\begin{theorem}\label{synchronous}  A group has a prefix-closed
  synchronous automatic structure with uniqueness if and only if it is
  isomorphic to the quotient of a finitely generated free group by a
  normal subgroup admitting a linear language of freely reduced generators with
  significant letters satisfying either of the
  following two conditions.
  \begin{enumerate}
  \item\label{k-central} For some constant $k$ each significant
  letters is $k$-central. 
  \item\label{o-central} The significant letters are $o$-central.
\end{enumerate}
\end{theorem}

A significant letter is $k$-central if it is within a distance $k$ of
the center of its word. A set of words with significant letters has
$o$-central significant letters if either the set is finite or if the
distance from significant letters to the center of their words is $o$
of the length of the word. Thus $k$-central implies $o$-central.

The generators mentioned in Theorems~\ref{asynchronous}
and~\ref{synchronous} are essentially the Schreier generators
corresponding to the combing from the automatic structure. The desired
properties of these generators are derived in a straightforward way
from the combing, but the argument in the other direction is more
complicated.

\section{Preliminary Definitions and Results} %%%%%%%%%%%%%%%%%%%%%%%%%%%

\subsection{Formal languages} %%%%%%%%%%%%%%%%%%%%%%%%%%%%%%%%%%%%%%%%

An alphabet is a finite nonempty set, $\Si$. A formal language over
$\Si$ is a subset of $\SSi$, the free monoid over $\Si$. Elements of
$\SSi$ are called words. The identity element of $\SSi$ is the empty
word, denoted $\e$. $\abs w$ is the length of a word $w$. If
$w=a_1\cdots a_n$, the distance from the letter $a_i$ to the center of
$w$ is $\abs{i - (n+1)/2}$. 

$\SSi$ is well ordered by the shortlex order, which is defined by
$u<v$ if either $\abs u < \abs v$ or $\abs u = \abs v$ and $u$ is less
than $v$ in the lexicographic order corresponding to some fixed
ordering of $\Si$. The shortlex order has the property that $u<v$
implies $xuy<xvy$ for all $x,y\in \SSi$.

We assume the reader is familiar with the theory of automatic groups
including the basic facts about regular languages and finite automata
necessary for the development of that theory. We include some
additional definitions and results from 
formal language theory which we shall need. See~\cite{RS} for a survey
of the whole field.

Recall that regular languages are the languages accepted by finite
automata. A finite automaton $\mathcal A$ over $\Si$ is a finite
directed graph with edge labels from $\Si_\e=\Si\cup\set{\e}$, a
designated initial vertex, and some terminal vertices. A path in
$\mathcal A$ is called successful if it starts at the initial vertex
and ends at a terminal vertex. The language accepted by $\mathcal A$
is the collection of labels of successful paths. The label of a path
is just the product of its edge labels. The label of a path of length
$0$ is $\e$. We assume without loss of generality that every edge and
every vertex of an automaton occur in some successful path. Other
edges and vertices can simply be deleted. We denote by $\abs\A$ the
number of vertices in a finite automaton $\A$.

Regular languages are closed under union, product and
generation of submonoid. We require some additional properties.

\begin{lemma}\label{regular} Every regular language $R$ may be
  expressed as a finite union $R=\cup X_iY_i$ of products of regular languages
  $X_i, Y_i$ in such a way that $w=uv\in R$ if and only if $u\in X_i,
  v\in Y_i$ for some $i$.
\end{lemma}

\begin{proof} Let $R$ be accepted by an automaton $\mathcal A$ with vertices
  $p_1,\ldots, p_n$. For each $i$ between $1$ and $n$ define an
  automaton $\mathcal A_i$ by altering $\mathcal A$ so that $p_i$ is
  its single terminal vertex. Likewise define $\mathcal A_i'$ by
  making $p_i$ the initial vertex. The languages $X_i$ accepted by
  $\mathcal A_i$ and $Y_i$ accepted by $\mathcal A_i'$ are as
  required.
\end{proof}

\begin{definition} For any word $w$, $w^r$ is $w$ written
  backwards. Likewise for any language $L$, $L^r=\set{w^r\mid w\in L}$
\end{definition}

\begin{lemma}\label{reverse} If $R$ is a regular language, so is $R^r$.
\end{lemma}

\begin{proof}
Let $R$ be accepted by the automaton $\mathcal A$. Reverse the
orientation of the edges of $\mathcal A$ and make the initial vertex
the single terminal. Taking each original terminal vertex in turn as
the initial vertex, we obtain a set of finite automata. The union of
the languages accepted by these automata is $R^r$.
\end{proof}

\subsection{Transductions}
Finite automata over $\Si\times\Si$ are defined just like finite
automata over $\Si$ except that edge labels are from $\Si_\e\times
\Si_\e$. The label of a path of length $0$ is $(\e,\e)$, and the
collection of labels of successful paths is a subset of $\SSi\times
\SSi$ called a rational transduction over $\Si$. A rational
transduction is a binary relation on $\SSi$.

\begin{lemma}\label{transduction} Every finite binary relation on $\SSi$
  is a rational transduction. Projections of rational transduction onto either
  coordinate yield regular languages. Rational transductions are closed under
  union and product. They are also closed under intersection with
  direct products $R\times S\subset \SSi \times \SSi$ of regular
  languages $R, S$ over $\Si$. 
\end{lemma}

\begin{proof}
The first two assertions are immediate from the definition of rational
transduction. 
To show closure under union combine two automata over $\Si\times\Si$
by adding a new initial vertex together with edges labeled $(\e,\e)$
from the new vertex to the initial vertex of each automaton. The new
automaton accepts the union of the two rational transductions accepted
by the original automata. Closure under product is demonstrated
similarly using edges from the terminal vertices of the first
automaton to the initial vertex of the second.

To complete the proof of the lemma it suffices to show that if $\rho$ is
a rational transduction accepted by an automaton $\A$ over $\Si\times
\Si$ and $R$ is a regular language accepted by the automaton $\B$ over
$\Si$, then $\rho\cap(\SSi\times R)$ and $\rho\cap(R \times \SSi)$ are
both rational transductions. The argument is the same in both cases.
We will show that $\rho\cap(\SSi\times R)$ is a rational
transduction. 

First observe that it does no harm to require a loop (an edge from a
vertex to itself) with label $(\e,\e)$ at each vertex of $\A$ and a
loop with label $\e$ at each vertex of $\B$. Now define an automaton
$\C=\A\times\B$ over $\Si\times \Si$ as follows. The set of vertices
of $\C$ is the Cartesian product of the vertices of $\A$ with the
vertices of $\B$. There is an edge with label $(a,b)$ from $(p_1,q_1)$
to $(p_2,q_2)$ if and only if there is an edge from $p_1$ to $p_2$
with label $(a,b)$ in $\A$ and an edge from $q_1$ to $q_2$ with label
$b$ in $\B$. 

It is easy to see that if there is path in $\C$ with label $(u,v)$
from from $(p_1,q_1)$ to $(p_2,q_2)$, then there is a path in $\A$
from $p_1$ to $p_2$ with label $(u,v)$ and a path in $\B$ from $q_1$
to $q_2$ with label $v$. The converse is also straightforward once we
observe that if there are paths in $\A$ from $p_1$ to $p_2$ with label
$(u,v)$ and in $\B$ from $q_1$ to $q_2$ with label $v$, then by
judiciously inserting loops with labels $(\e,\e)$ or $\e$ into the two
paths we can arrange things so that $v$ is expressed in exactly
the same way as a product of elements of $\Si_\e$ along
both paths. (This argument is given in greater detail in the proof
of~\cite[Theorem 4.4]{Gi}.)

Take the initial vertex of $\C$ to be $(p_0,q_0)$ where $p_0$ is the
initial vertex of $\A$ and $q_0$ is the initial vertex of
$\B$. Likewise $(p,q)$ is terminal if both $p$ and $q$ are.  It
follows from the preceding paragraph that $(u,v)$ is the label of a
successful path in $\C$ if and only if $(u,v)$ is the label of a
successful path in $\C$ and $v$ is the label of a successful path in
$\B$.
\end{proof}

\begin{lemma}\label{identity} For each regular language $R$ over $\Si$
  the binary relation $\rho_R=\set{(u,u)\mid u\in R}$ is a rational
  transduction.
\end{lemma}
\begin{proof} Since $\rho_R=\rho_\SSi\cap (R\times R)$, it suffices by
  Lemma~\ref{transduction} to consider the case $R=\SSi$. The
  automaton with one vertex $p$ (which is both initial and terminal)
  and edges $p \edge{(a,a)} p$ for each $a\in \Si$ accepts $\rho_\SSi$
\end{proof}

\subsection{Linear languages} %%%%%%%%%%%%%%%%%%%%%%%%%%%%%%%%%%%%%%

\begin{definition}\label{linear-definition} A language $L$ over $\Si$
  is linear if for some rational transduction $\rho$ over $\Si$,
  $L=\set{uv^r\mid (u,v)\in \rho}$.
\end{definition}

In other words a linear language consists of all words $uv^r$ such
that $(u,v)$ is the label of a successful path in some fixed automaton
over $\Si\times\Si$. Other characterizations are given
in~\cite[Chapter 3, Section 6.1]{RS}. Automata over
$\Si\times\Si$ serve as acceptors for both rational transductions and
linear languages.

\begin{lemma}\label{linear-properties} The union of two linear
  languages is linear.  The intersection of a linear
  language over $\Si$ and a regular language over $\Si$ is linear. 
\end{lemma}

\begin{proof} The first assertion is immediate from
  Lemma~\ref{transduction}. For 
  the second let $L$ be linear and $R$ regular. $L$ is accepted by an
  automaton $\mathcal A$ which also accepts a rational transduction
  $\rho$ such that $L = \set{uv^r \mid (u,v)\in\rho}$. Express $R=\cup
  X_iY_i$ as in 
  Lemma~\ref{regular}. By Lemmas~\ref{reverse} and~\ref{transduction}
  $\rho'=\cup (\rho \cap (X_i\times Y_i^r))$ is a rational
  transduction. Since $L\cap R=\set{uv^r\mid (u,v)\in \rho'}$,
  $L\cap R$ is linear.
\end{proof}  

\subsection{Languages and groups} %%%%%%%%%%%%%%%%%%%%%%%%%%%%%%%%%%%%

Consider a group $G$ and a surjective homomorphism $\mu:F\to G$ from a
finitely generated free group $F$. Let $N$ be the kernel of
$\mu$. Take $\Si$ to be an alphabet of free generators and their
inverses for $F$ and let $\pi:\SSi \to F$ be the projection which
sends each word $w\in \SSi$ to the element of $F$ it
represents. Notice that $\SSi$ is equipped with formal inverses in a
natural way, and $\pi$ respects inverses. We will call this
configuration a choice of generators for $G$. From now on $\Si$ stands
for an alphabet with formal inverses.
\begin{equation}\label{cog}
\SSi\edge{\pi}F\edge{\mu}G.
\end{equation}
When we wish to avoid explicit reference to $\mu$ and $\pi$, we will
use $\hatt x$ and $\ovr x$ to denote the image of $x$ in $F$ and $G$
respectively.

Given a choice of generators~(\ref{cog}), we see that for every
language $L$ over $\Si$ there is a subgroup $H=\subgroup{\hatt L}$
generated by the image of $L$ in $F$. We call $L$ a language of
generators for $H$.

\subsection{Significant letters}

\begin{definition}\label{significant} Let $L\subset\SSi$ be a
  language of freely reduced words which does not contain the empty
  word.  $L$ has significant letters if every $w \in L$ can be written
  as a product $w=uav\inv$ with $a\in \Si$ such that for all
  $w_1,w_2\in L$ and $\epsilon_1,\epsilon_2=\pm 1$, free reduction of
  $(w_1)^{\epsilon_1}(w_2)^{\epsilon_2} =
  (u_1a_1v_1\inv)^{\epsilon_1}(u_2a_2v_2\inv)^{\epsilon_2}$ does not
  affect $a_1$ or $a_2$ unless the product reduces to $\e$.
\end{definition}

When considering a word $w$ in a language $L$ with significant
letters, $w=uav\inv$ will always mean the significant letter
decomposition of $w$. Significant letters need not be uniquely
determined, but it is clear from Definition~\ref{significant} that we
may assume that if $w,w\inv\in L$ and $w=uav\inv$, then $w\inv=va\inv
u\inv$ is the significant letter decomposition of $w$. It follows that
if $L$ has significant letters, so does $L\cup L\inv$. We record this
fact along with two immediate consequences of Definition~\ref{significant}.

\begin{lemma}\label{overlap} Let $L$ have significant letters. 
  Then $L\cup L\inv$ has significant letters.
  Consider $w=uav\inv, w_1\in L$. If either $ua$ is
  a prefix of $w_1$ or $av\inv$ is a suffix, then $w=w_1$. If either
  $va\inv$ is a prefix or $a\inv u$ is a suffix, then $w=w_1\inv$.
\end{lemma}

If $L \subset \SSi$ has significant letters, then $\hatt L$ is a set
of free generators for the subgroup
$\subgroup{\hatt L}\subset F$ generated by $\hatt L$. Indeed if
$w_1^{\epsilon_1}\cdots w_n^{\epsilon_n}$ is any product of words from
$L$ and their inverses such that no $w_i^{\epsilon_i}$ is
followed by its inverse, then free reduction of the product does not
affect the significant letter of any $w_i$. Consequently
$\hatt{w_1}^{\epsilon_1}\cdots \hatt{w_n}^{\epsilon_1}\ne 1$.

\subsection{Combings}

\begin{definition}\label{combing} A combing is a language $C$ over
  $\Si$ such that $\ovr C =G$.  $C$ is prefix closed if every prefix
  of any $w\in C$ is also in $C$. $C$ is a combing with uniqueness if
  $C$ maps bijectively to $G$. $C$ is regular if it is a regular
  language.
\end{definition}

There are other definitions of combing in the literature.

\begin{lemma}\label{subwords} If $C$ is a prefix closed combing with
  uniqueness, then no nontrivial subword of a word in $C$ defines the
  identity in $G$. In particular $C$ consists of freely reduced words.
\end{lemma}

\begin{proof} If not, then there is $uxv\in C$ with $\ovr
  x = 1$. By closure under prefixes, $u, ux \in C$, contradicting
  uniqueness.
\end{proof}

\begin{lemma}\label{linear-inverse} If $\rho$ is a rational
  transduction over $\Si$, then $L = \set{uv\inv\mid (u,v)\in \rho}$
  is a linear language. If $L$ is linear, then so is $L\inv=\set{w
  \mid w\inv\in L}$.
\end{lemma}
\begin{proof} For the first part replace each edge label $(ab)$ with
  $(a,b\inv)$.  For the second assertion change each label $(a,b)$ to
  $(b\inv ,a\inv)$.
\end{proof} 

In practice we will not replace $(a,b)$ by $(a,b\inv)$. Instead we
will just read $(uv\inv)$ instead of $(uv^r)$ for each path with label
$(u,v)$ in an automaton accepting $\rho$. From now on the linear
language corresponding to a transduction $\tau$ will be
$L=\set{(uv\inv\mid (u,v) \in \tau}$.

We will make use of the following possibly infinite automaton.

\begin{definition}\label{Cayley} Let $G$ be a group and~(\ref{cog}) a
  choice of generators. The Cayley automaton $\mathcal A_G$ has vertices
  $G$ and edges $g \edge{(a,b)} h$ for all $g,h\in 
  G$ and $a,b\in \Si_\e$ with $g \ovr b = \ovr a h$. The initial state
  of $\mathcal A_G$ is $1$, and all states are terminal states. 
\end{definition}

\begin{lemma}\label{ball} There is a path in $\A_G$ 
  with label $(u,v)$ from $1$ to $h$ if and only if $h=\ovr {u\inv
  v}$. If $u$ and $v$ are asynchronous $k$-fellow travelers, then
  the path may be chosen in the ball of radius $k$ around $1$.
\end{lemma}

\begin{proof}
It is straightforward to prove by induction on length that there is a
path with label $(u,v)$ from $g$ to $h$ in $\A_G$ if and only if
$g\ovr u= \ovr v h$. The second assertion follows directly from the
definition of asynchronous fellow traveler.
\end{proof}

\subsection{Automatic structures} %%%%%%%%%%%%%%%%%%%%%%%%%%%%%%%%%%%

\begin{lemma}\label{asynchronous-definition} A combing $C$ for a group
  $G$ supports a prefix-closed asynchronous automatic structure with
  uniqueness if and only if $C$ is prefix-closed with uniqueness, and for each
  $a\in \Si$ the binary relation $\rho_a=\set{(u,v) \mid u,v\in C,
  \ovr u\ovr a=\ovr v}$ is a rational transduction.
\end{lemma}

\begin{proof} Suppose $C$ supports a prefix-closed asynchronous automatic
  structure with uniqueness in the sense of~\cite[Definition
  7.2.1]{E+}. Then $C$ is a prefix-closed combing with uniqueness, and
  it is not hard to check that the corresponding binary relations
  are transductions.

  For the converse take $C$ and $\rho_a$, $a\in \Si$ as above, and
  let $k$ be an upper bound for the number of vertices in automata
  $\mathcal A_a$ accepting $\rho_a$. Suppose $\gamma$ is a
  successful path in some $\mathcal A_a$. For each vertex $p$ of
  $\gamma$ there is a path of length at most $k$ from $p$ to a
  terminal vertex of $\mathcal A_a$. Thus if $(u,v)$ is the label of
  $\gamma$ up to $p$, there are words $x,y$ of length at most $k$ such
  that $(ux,vy)\in \rho_a$. Consequently $\ovr{uxa}=\ovr{vy}$, which
  implies that the word difference $u\inv v$ has the same image in $G$
  as some word of length at most $2k+1$. From this observation
  together with the fact that $\rho_\e$ is the identity binary relation
  on $C$ we see that $C$ satisfies the asynchronous fellow traveler
  property. By Theorems 1 and 2 of~\cite{MS} some subset of $C$ is a
  regular combing supporting an asynchronous automatic structure.
  Since $C$ is a combing with uniqueness, the subset must be $C$ itself.
\end{proof}

An analog of Lemma~\ref{asynchronous-definition} holds for synchronous
automatic structures, but the rational transductions are of a special
type.

\begin{definition}\label{synchronized} A rational transduction
  $\rho\subset \SSi\times 
  \SSi$ is called synchronized if $(u,v)\in \rho$ implies the lengths
  $\abs u$ and $\abs v$ differ by at most $k$ for some constant $k$. A
  finite automaton over $\Si\times \Si$ is synchronized if it is built
  up from a subautomaton $\A_0$ with edge labels all in
  $\Si\times\Si$ by attaching directed paths of length at most $k$ such
  that the edge labels along each path are either all in $\Si\times
  \set\e$ or all in $\set\e\times\Si$. These paths are attached at
  their initial points only and are otherwise disjoint from each
  other. 
\end{definition}

It is clear that any rational transduction accepted by a synchronized
automaton is synchronized. The converse follows from
\cite[Proposition~2.1]{FS}.  

\begin{lemma}\label{synchronized-automaton} A rational transduction is
  synchronized if and only if it is accepted by a synchronized finite
  automaton. 
\end{lemma}

\begin{lemma}\label{synchronous-definition} A combing $C$ for a group
  $G$ supports a prefix-closed synchronous automatic structure with
  uniqueness if and only if $C$ is prefix-closed with uniqueness and for each
  $a\in \Si$ the binary relation $\rho_a=\set{(u,v) \mid u,v\in C,
  \ovr u\ovr a=\ovr v}$ is a synchronized rational transduction.
\end{lemma}

\begin{proof} Suppose $C$ supports a prefix-closed synchronous automatic
  structure with uniqueness in the sense of~\cite[Definition
  2.3.1]{E+}. It is clear that that associated binary relations
  $\rho_a$ are rational transductions. By~\cite[Lemma 2.3.9]{E+} the
  uniqueness condition on $C$ implies that the $\rho_a$'s are
  synchronized rational transductions.

  For the converse take $C$ and $\rho_a$, $a\in \Si$ as above.
  Synchronized finite automata accepting the $\rho_a$'s fit the
  definition of the the automata occurring in ~\cite[Definition
  2.3.1]{E+} once labels $(a,\e)$ and $(\e, a)$ are replaced by labels
  $(a,\$)$ and $(\$, a)$ respectively. The same conclusion holds for
  $\rho_\e$, as it is the identity on $C$.
\end{proof} 

Automatic structures can also be defined in terms of regular combings
satisfying fellow traveler conditions. These conditions are defined in
terms of the word metric $d$ corresponding to a choice of
generators~(\ref{cog}). We write
$D_a(w,v)\le k$ if two words $w,v\in \SSi$ satisfy the asynchronous $k$-fellow
traveler condition and $D_s(w,v)\le k$ if they satisfy the
synchronous $k$-fellow traveler condition. The following lemma
records some well known properties.

\begin{lemma}\label{fellow-traveller} If $D_s(w,v)\le k$, then
  $D_a(w,v)\le k$. Further
\begin{enumerate} 
\item If $D_a(u,v)\le k$ and $D_a(v,v')\le k'$, then
  $D_a(u,v')\le k+k'$.  
\item If $D_s(u,v)\le k$ and $D_s(v,v')\le k'$, then
  $D_s(u,v')\le k+k'$.  
\item $D_s(u, uv)\le \abs{v}$.
\end{enumerate}
\end{lemma}

\section{Finding Generators} %%%%%%%%%%%%%%%%%%%%%%%%%%%%%%%%%%%%%%%%%

In this section we prove Theorems~\ref{asynchronous}
and~\ref{synchronous} in one direction by extracting from an automatic
structure a language of generators of the required type. The arguments
are identical for both types of automatic group except for one
paragraph which applies only to the synchronous case.

Let $G$ be automatic of either type. Make a choice of
generators~(\ref{cog}), and take $N$ to be the kernel of $\mu$. As in
Definitions~\ref{asynchronous-definition}
and~\ref{synchronous-definition} $C$ is a combing supporting a prefix
closed automatic structure structure with uniqueness, and for each
$a\in \Si$, $\rho_a=\set{(u,v) \mid u,v\in C, a\in \Si, \ovr u \ovr a
= \ovr v}$ is a rational transduction. In the synchronous case
$\rho_a$ is a synchronized rational transduction. We will show that
$L=\set{ uav\inv \mid u,v\in C, a\in \Si, \ovr{ua} = \ovr v, uav\inv
\text{ is freely reduced}}$ is the desired language of generators. 

First we note that by construction $L$ is closed under inverse. Next
we show that $L$ is a linear language. By
Lemma~\ref{transduction} the product $(\rho_a)\set{(a,\e)}=\set{(ua,v)
\mid u,v\in C, a\in \Si, \ovr u \ovr a = \ovr v}$ is a rational
transduction. Likewise $\rho=\cup_{a\in\Si}\rho_a=\set{(ua,v) \mid
  u,v\in C, a\in \Si, \ovr 
u \ovr a = \ovr v}$ is also a rational transduction. Hence
$L'=\set{uav\inv\mid u,v\in C, a\in \Si, \ovr u \ovr a = \ovr v}$ is a
linear language. As $L$ is the intersection of $L'$ with the regular
language of nontrivial freely reduced words, $L$ is linear by
Lemma~\ref{linear-properties}.

In the synchronous case $\rho$ is synchronous because each $\rho_a$
is. Thus for some positive integer $k$, $(ua,v)\in \rho$ implies that
$\abs{ua}$ and $\abs v$ differ by at most $k$. We conclude that in the
synchronous case the $a$'s are $k$-central for words in $L$ and hence
$o$-central as well.

It remains to show that in both cases $L$ is a language of generators
and the $a$'s are significant letters for $L$. Observe that prefix
closure and uniqueness for $C$ imply that $\hatt C$ is a set of prefix
closed coset representatives for $N$ in $F$. We will interpret this fact
geometrically. 

Each $w\in \SSi$ may be thought of as a path beginning at $1$ in $\G$,
the Cayley diagram of $G$ with respect to the set of generators
$\Si$. We pick one letter from each pair $a,a\inv$ to use as edge
labels in $\G$. An edge of $\G$ traversed backwards is construed as a
forward edge with the inverse label. $C$ is a spanning tree for $\G$, and
any word $uav\inv$ with $\ovr u \ovr a = \ovr v$ is a cycle. If $a$
labels an edge of $\G$ in the spanning tree $C$, then because of our
convention about edge labels, $uav\inv$ is a cycle
in $C$ and thus freely equal to $\e$. Otherwise $uav\inv$ is freely
reduced by inspection. Likewise free reduction of a product of two
words in $u_1a_1v_1\inv, u_2a_2v_2\inv\in L$, cannot involve the
$a_i$'s unless they are labels of inverse edges in $\G$ in which case
the product is freely equal to $\e$. Finally if $uav\inv$ is freely
reduced, so is $va\inv u\inv$.  It follows that the $a$'s are
significant letters for $L$.  Hence $L$ is a language of free
generators and their inverses for the subgroup $\subgroup{\hatt L}\subset F$.

A word $w\in \SSi$ represents an element of $N$ if and only if $w$ is
a cycle in $\G$. Thus $\subgroup{\hatt L}\subset N$. On the other hand
suppose $w$ is a cycle in $\G$. A short argument by induction on the number,
$n$, of edges of $w$ not in $\G$ shows that $w$ is freely equal to a
product of words in $L$. Indeed if $n=0$, then as above $w$ is a cycle
in $C$ and so freely equal to $\e$, which is the empty
product.  Otherwise $w=uax$ where $a$ labels the first edge not in
$C$. But then $u \in C$, and there is $uav\inv\in L$. Consequently $w$
is freely equal to $(uav\inv)vx$. But $vx$ is a cycle to which the
induction hypothesis applies.

\section{Finding Automatic Structures} %%%%%%%%%%%%%

We complete the proofs of Theorems~\ref{asynchronous}
and~\ref{synchronous} by finding the required automatic structures.
Assume that $G$ is a group with choice of generators~(\ref{cog}) and
that $N$ has a linear language, $L$, of freely reduced generators with
significant letters. By Lemmas~\ref{overlap} and~\ref{linear-properties}
we may assume that $L$ is closed under inverse.

It suffices to show that $\SSi$ contains a prefix closed regular
combing with uniqueness which
satisfies the appropriate fellow-traveler property. The arguments in
the two cases are almost identical.  When it is necessary to
distinguish between them, we refer to the central and non-central
cases.

Let $L$ be accepted by an automaton $\A$ over $\Si_\e\times \Si_\e$.
If possible choose $\A$ to be synchronized. Let $\tau$ be the rational
transduction accepted by $\A$; $L=\set{ uv\inv \mid (u,v)\in \tau}$. 

If $\A$ is synchronized, there are no edges with label
$(\e,\e)$. However, in general there may be some. If there is a cycle
with label $(\e,\e)$, then identifying all the vertices in the cycle,
discarding the edges in the cycle, and taking the resulting vertex to
be initial or terminal if one of the identified vertices was does not
change the set of labels of successful paths. Consequently we assume
there are no such cycles. 

Without loss of generality delete edges and
vertices of $\A$ not lying on successful paths. If $\A$ was
synchronized before this change, it remains so. Choose a constant $K$
greater than the number of vertices and edges in $\A$.

If $G$ is free or finite, there is nothing to prove. Thus we may
assume $N\ne 1$ and $N$ has infinite index in $F$. As finitely
generated normal subgroups of free groups have finite index, $N$ is
not finitely generated. Thus $L$ is infinite, and consequently $\A$
has at least one cycle.

Define $\A_0$ to be the subgraph of $\A$ consisting of all vertices 
and edges which are in cycles or in paths leading to cycles.
As every edge of
$\A$ is on a successful path, the initial vertex of $\A$ must be in
$\A_0$. By definition of $\A_0$ there are no cycles in $\A-\A_0$ and no
edges from $\A-\A_0$ into $\A_0$. Consequently every
path in $\A$ lies in $\A_0$ except for its last $j$ vertices for some
$j\le K$. 

\begin{lemma}\label{upto}
An edge of $\A$ whose label holds the significant letter for some
successful path does not lie in $\A_0$. 
\end{lemma}

\begin{proof}
Assume otherwise. There is a successful path $\gamma=\gamma_1\gamma_2$
such that $\gamma_1$ is in $\mathcal A_0$ and the significant letter
occurs in an edge label of $\gamma_1$.  Since $\gamma_1$ is in
$\mathcal A_0$, there are successful paths
$\gamma_1\gamma_3\gamma_4^i\gamma_5$ where $\gamma_4$ is a cycle. By
Lemma 2.9 all these paths have the same label.  But then $\gamma_4$
must have label $(\e,\e)$ contrary to our assumption about cycles in
$\A$.
\end{proof}

\begin{lemma}\label{synchronous-transduction}
In the central case $\tau$ is a synchronized rational transduction, and
the edge labels of $\A_0$ lie in $\Si\times \Si$. 
\end{lemma}

\begin{proof} 
We claim that
every cycle in $\A$ has label $(u,v)$ with $\abs u = \abs v$. As every
path in $\A$ has at most $K$ edges which do not occur in cycles along
the path, it will follow that $\abs{\abs x - \abs y}\le K$ for every
$(x,y)\in \tau$. Hence $\tau$ will be synchronized. 

To verify our claim suppose $(u,v)$ is the label of a cycle in $\A$
and $\abs u\ne \abs v$. As all cycles lie in $\A_0$, Lemma~\ref{upto}
implies that for fixed words $u_0,u_1, v_0, v_1$ and all integers $i\ge
0$, $L$ contains words $u_0u^iu_1v_1\inv v^{-i}v_0\inv$ whose significant
letters occur in the subword $u_1v_1\inv$. It follows by a
straightforward argument that the significant letters of $L$ are not
$o$-central and hence not $k$-central. Thus our claim is valid.

Finally since $\tau$ is synchronized, our choice of $\A$ guarantees
that $\A$ is too. It follows from Definition~\ref{synchronized} that
the edge labels of $\A_0$ lie in $\Si\times \Si$.  
\end{proof}

The choice of generators~(\ref{cog}) determines a Cayley diagram
$\Gamma$ for $G$ with the corresponding word metric $d$. Each word
in $\SSi$ is the label of a unique path from $1$ in $\Gamma$, and we
will use $w$ to refer to the path as well as the word. A word
represents an element of $N$ if and only if it is a cycle in $\Gamma$.

\begin{lemma}\label{first-combing} The language $C$ consisting of all
  prefixes not including the significant letter of each
  $w\in L$ is a prefix closed combing with uniqueness for $G$. Further
  if $\ovr{ua}=\ovr v$ for $u,v\in C$ and $a\in \Si$, then either
  $uav\inv$ is freely equal to $\e$ or $uav\inv \in L$ with
  significant letter $a$.
\end{lemma}

\begin{proof} $C$ is obviously prefix closed.
  For any $g\in G$ there is a simple path $w$ in $\G$ from $1$ to
  $g$. Since $N\ne 1$, there is a simple cycle of length at least $1$
  starting at $g$. Extend the path $w$ by continuing around this cycle
  until its first return to $w$ and then following $w$ back to
  $1$. This extension of $w$ is a cycle passing through $g$ with
  freely reduced label. Consequently $w$ is the free reduction of a product of
  generators from $L$. Each of these generators is a cycle, and one of
  them, say $w=uav\inv$, must contain the vertex $g$. Consequently
  some prefix of $u$ or of $v$ is a path from $1$ to $g$. As $L$ is
  closed under taking inverses, that prefix lies in $C$. Thus $C$ maps
  onto $G$.

  To prove that $C$ maps injectively to $G$ suppose $\ovr u =
  \ovr v$ for $u,v\in C$ with $u\ne v$. It follows that 
  $uv\inv $ is freely equal to a nonempty product of generators from
  $L$. By the nature of significant letters, the prefix $u_1a_1$ from
  the first generator $u_1a_1v_1\inv$ in the product and the suffix
  $a_nv_n\inv$ from the last generator are not affected by free
  reduction of the product. As $u$ and $v$ are both freely reduced, it
  follows that $u_1a_1$ is a prefix of $u$ or $a_nv_n\inv$ is a suffix
  of $v\inv$. But by the definition of $C$ together with
  Lemma~\ref{overlap} this is impossible.
   
  The last assertion is proved in the same manner. As $u$ and $v$ are
  freely reduced, $uav\inv$ is either freely reduced or freely equal
  to $u_1v_1\inv$ for prefixes $u_1$ of $u$ and $v_1$ of $v$. In the
  latter case $u_1=v_1$ by injectivity whence $uav\inv$ is freely
  equal to $\e$. In the former case the argument of the previous
  paragraph yields either $ua=u_1a_1$ or $av\inv=a_nv_n\inv$. It
  follows that $n=1$ and $uav\inv=u_1a_1v_1\inv$.
\end{proof}

\begin{lemma}\label{ft}  In the
  non-central case there is a constant $k$ such that $C$ satisfies the
  asynchronous $k$-fellow traveler condition. In the central case
  there is a constant $k$ such that $C$ satisfies the synchronous
  $k$-fellow traveler condition.
\end{lemma}

\begin{proof} Suppose $u,v\in C$ with $d(\ovr u, \ovr v)\le 1$ in
  $\G$. If $d(\ovr u, \ovr v)=0$, then $u=v$ by uniqueness and both
  fellow traveler conditions are satisfied with $k=0$. If $d(\ovr u,
  \ovr v)=1$, then $\ovr{ua}=\ovr v$ for some $a\in \Si$. By
  Lemma~\ref{first-combing} either $ua$ is freely equal to $v$ or
  $uav\inv \in L$. In the first case both fellow traveler conditions
  are satisfied with $k=1$.
 
  Assume the second case holds, and suppose $\gamma$ is a successful
  path in $\mathcal A$ with label $uav\inv$. $\gamma$ consists of a prefix
  $\gamma_0$ in $\A_0$ followed by a suffix of length at most
  $K$. Let $(u_0,v_0)$ be the label of $\gamma_0$. By
  Lemma~\ref{upto} $u_0$ includes all but the last
  $j$ letters of $u$ for some $j\le 2K$, and likewise for
  $v_0$. In the central case $\abs{u_0}=\abs{v_0}$ as the edge labels
  of $\A_0$ are all from $\Si\times\Si$. By
  Lemma~\ref{fellow-traveller} it suffices to prove 
  that $u_0$ and $v_0$ are $k_0$ fellow travelers of the appropriate
  type for some constant $k_0$.

  Consider any vertex $p$ of $\gamma_0$. There is a path of length at
  most $K$ from $p$ to a terminal vertex of $\A$. Thus
  $u_0x(v_0y)\inv \in L$ for some words $x,y$ with $\abs x, \abs y \le
  K$. Consequently $\ovr{u_0x}=\ovr{v_0y}$, which implies that
  the word difference $u_0\inv v_0$ has the same image in $G$ as some
  word of length at most $2K$. In the non-central case we see
  immediately that $D_a(u_0,v_0)\le 2K$. In the central case
  $D_s(u_0,v_0)\le 2K$ because the edge labels of $\A_0$ are all
  in $\Si\times\Si$.
\end{proof}

It remains only to show that $C$ is regular, but unfortunately
there does not seem to be any reason why this should be so. However,
by replacing certain suffixes of length at most $2K$ of words in
$C$ with new suffixes of length at most $2K$ we
obtain a combing $C'$ which works. 

Recall that $\A_0$ is the subgraph of $\A$ supported by all vertices
which are in cycles or in paths leading to cycles and that $\A_0$
contains the initial vertex of $\A$. Make $\A_0$ into an automaton
$\B_0$ over $\Si$ by replacing each edge label $(a,b)$ with the label
$a$. The initial vertex of $\B_0$ is the initial vertex of $\A$, and
all vertices are terminal. $\B_0$ accepts a prefix closed regular
language $C_0$. By Lemma~\ref{upto} $C_0$ is a collection of prefixes
of $C$. It follows from the structure of $\A$ that each word in $C$ is
obtained by appending a word of length at most $2K$ to a word in
$C_0$. We will define $C'$ by appending other suffixes of at most the
same length.

Let $X$ be the set of all words in $\SSi$ of length at most $2K$.
Clearly $C\subset C_0X$. Define $C'$ as follows.  For each $g\in G$
pick the unique $x\in X$ minimum in the shortlex order such that there
exists $u_0\in C_0$ with $\ovr{u_0x}=g$. Since $C\subset C_0X$, such a
$u_0$ exists. By the uniqueness property of $C$, there is just one
choice for $u_0$.  Also since $\e$ is the minimum element of $X$ in
the shortlex order, our construction guarantees $C_0\subset C'$.

\begin{lemma}\label{second-combing}
$C'$ is a prefix closed combing with uniqueness. For some constant
  $k'$, $C'$ satisfies the appropriate $k'$-fellow traveler
  condition.
\end{lemma}

\begin{proof}
  $C'$ has uniqueness by construction. Likewise Lemma~\ref{ft} and the
  properties listed in Lemma~\ref{fellow-traveller} insure that $C'$
  satisfies the appropriate fellow traveler condition. To show prefix
  closure consider a prefix $v$ of $u_0x\in C'$. If $v$ is a prefix of
  $u_0$, then $v\in C_0\subset C'$. Otherwise $v = u_0x_1$ for some
  prefix $x_1$ of $x=x_1x_2$. If $v\notin C'$, then there exists
  $u_1y\in C'$ with $u_1\in C_0$, $\ovr{u_1y}=\ovr{u_0x_1}$ and $y <
  x_1$. But then $yx_2< x$ and $\ovr{u_1yx_2}=\ovr{u_0x}$
  contradicting the construction of $C'$.
\end{proof}

We must show that $C'$ is a regular language. 
For each $x\in X$ let $C_x=\set{r\mid r\in C_0, rx\in C'}$. $C'=\cup_x
  C_xx$ is regular if each $C_x$ is. $C_x=C_0-\cup_{y\in X, y<x}
  C_{x,y}$ where $C_{x,y}= \set{r\mid r\in
    C_0, \ovr{rx}=\ovr{sy} \mbox{ for some $s\in C_0$}}$, so it
  suffices to show that $C_{x,y}$ is regular.

Define a finite automaton over $\Si_\epsilon\times\Si_\epsilon$ from the
ball of radius $k+4K$ around $1$ in the Cayley automaton $\A_G$ by
taking $1$ as the initial vertex and $\ovr{xy\inv}$ as the single terminal
vertex. Let $\tau_{x,y}$ be the rational transduction accepted by this
automaton. By Lemma~\ref{ball} $\tau_{x,y}$ is contained in the set of
$(u,v)$ such that $\ovr {u\inv v} = \ovr{xy\inv}$ and contains all
$(u,v)$ such that$\ovr {u\inv v} = \ovr{xy\inv}$, and $D_a(u, v)\le k+4K$.

Suppose $r,s\in C_0$ with $\ovr{rx}=\ovr{sy}$ for $x,y\in
X$. As $C_0\subset C$, Lemmas~\ref{ft} and~\ref{fellow-traveller}
imply $D_a(r,s)\le k+\abs x + \abs y \le k+4K$. 
Hence $\tau_{x,y}\cap (C_0\times C_0)$ is a rational transduction
whose projection onto the first coordinate is $C_{x,y}$. Thus
$C_{x,y}$ is regular.

\end{document}